\documentclass{article}
\usepackage{array}
\usepackage{amsmath,amssymb,amsthm}
\numberwithin{equation}{section}

\def\R{\mathbb{R}}
\def\C{\mathbb{C}}

\def\muhat{{\widehat\mu}}
\def\Sigmahat{{\widehat\Sigma}}
\def\wtil{\widetilde}

\theoremstyle{plain}
\newtheorem{thm}{Theorem}[section]
\newtheorem{lemma}[thm]{Lemma}
\newtheorem{prop}[thm]{Proposition}

\newtheorem{rmk}[thm]{Remark}

\begin{document}
\renewcommand{\baselinestretch}{1.2}
\markright{
}
\markboth{\hfill{\footnotesize\rm Max-Louis G. Buot, Serkan 
Ho\c{s}ten, and Donald St. P. Richards
}\hfill}
{\hfill {\footnotesize\rm THE LIKELIHOOD EQUATIONS FOR THE 
BEHRENS-FISHER PROBLEM} \hfill}
\renewcommand{\thefootnote}{}
$\ $\par
\fontsize{10.95}{14pt plus.8pt minus .6pt}\selectfont
\vspace{0.8pc}
\centerline{\large\bf COUNTING AND LOCATING THE SOLUTIONS OF}
\vspace{2pt}
\centerline{\large\bf POLYNOMIAL SYSTEMS OF MAXIMUM LIKELIHOOD}
\vspace{2pt}
\centerline{\large\bf EQUATIONS, II: THE BEHRENS-FISHER PROBLEM}
\vspace{.4cm}
\centerline{Max-Louis G. Buot, Serkan Ho\c{s}ten, and 
Donald St. P. Richards} 
\vspace{.4cm} 
\centerline{\it Xavier University, San Francisco State University,}
\vspace{.2cm}
\centerline{\it Penn State University and SAMSI}
\vspace{.55cm}
\fontsize{9}{11.5pt plus.8pt minus .6pt}\selectfont

\begin{quotation}
\noindent {\it Abstract:}
Let $\mu$ be a $p$-dimensional vector, and let $\Sigma_1$ and 
$\Sigma_2$ be $p \times p$ positive definite covariance matrices.  
On being given random samples of sizes $N_1$ and $N_2$ from 
independent multivariate normal populations $N_p(\mu,\Sigma_1)$ 
and $N_p(\mu,\Sigma_2)$, respectively, the Behrens-Fisher problem 
is to solve the likelihood equations for estimating the unknown 
parameters $\mu$, $\Sigma_1$, and $\Sigma_2$.  We shall prove 
that for $N_1, N_2 > p$ 
there are, almost surely, exactly $2p+1$ complex solutions 
of the likelihood equations. 
For the case in which $p = 2$, we utilize Monte Carlo 
simulation to estimate the relative frequency with which a 
typical Behrens-Fisher problem has multiple real solutions; 
we find that multiple real solutions occur infrequently.
\par

\vspace{9pt}
\noindent {\it Key words and phrases:}
Behrens-Fisher problem, B\'ezout's theorem, maximum likelihood 
estimation, maximum likelihood degree.
\par
\end{quotation}\par

\fontsize{10.95}{14pt plus.8pt minus .6pt}\selectfont

\section{Introduction}
\label{sec:intro}
\setcounter{equation}{0}

Let $\mu \in \R^p$ be a $p$-dimensional vector, and let 
$\Sigma_1$ and $\Sigma_2$ be $p \times p$ positive definite 
(symmetric) matrices.  Consider independent multivariate 
normal populations, $N_p(\mu,\Sigma_1)$ and $N_p(\mu,\Sigma_2)$, from 
which we have been given random samples $X_1,\ldots,X_{N_1}$ and 
$Y_1,\ldots,Y_{N_2}$, respectively.   On the basis of the given data, 
the famous Behrens-Fisher problem (Behrens (1929), Fisher (1939)) is 
to estimate the parameters $\mu$, $\Sigma_1$, and $\Sigma_2$ by means 
of the method of maximum likelihood.  

It is well-known that the corresponding system of likelihood 
equations cannot be solved explicitly, and that has led many to 
propose alternative solutions to the Behrens-Fisher problem 
(Anderson, 2003, p. 187 ff.).  More importantly, the Behrens-Fisher 
problem is an early example of a 
hypothesis testing problem involving exponential families of 
densities and for which the resulting sufficient statistics, when the 
parameters are restricted to the parameter space determined by $H_0$, 
fail to be complete (Linnik, 1967).  In such a situation, nuisance 
parameters exist, and the construction of an exact size-$\alpha$ test 
is a difficult problem.  

Consequently, the literature on the 
Behrens-Fisher problem is substantial, reflecting the intense interest 
which the problem has generated since its inception.  Indeed, the 
problem has generated an extensive philosophical discussion as well 
as many efforts to derive solutions which are optimal for 
statistical inference (Wallace (1980), Kim and Cohen (1998), 
Stuart and Ord (1994)).  In this paper, we determine the number of 
solutions of the likelihood equations.  

For the case in which $p=1$, there are three unknown scalar 
parameters, {\it viz.}, $\mu$, 
the common mean, and $\sigma_1^2$ and $\sigma_2^2$, the population 
variances.  In this case, Sugiura and Gupta (1987) reduced the 
system of equations to a cubic equation in $\mu$ and deduced that, 
almost surely, there are three complex solutions; they observed also 
that the likelihood equation tended to have multiple real solutions 
if $\sigma_1^2$ and $\sigma_2^2$ are small in comparison with $\mu$, 
and otherwise that the likelihood equation usually has a unique real 
solution. Drton (2007) also studied the univariate Behrens-Fisher 
problem and showed, in particular, that if the null hypothesis is 
true then the probability of multiple real solutions tends to zero 
as the sample sizes tend to infinity. We will prove the analogous 
result for the multivariate Behrens-Fisher problem in Theorem 
\ref{thm:uniqueness}.

In this paper, as in the article of Buot and Richards (2006),  
we apply results from the 
theory of algebraic geometry to study the solution set of the system 
of likelihood equations for the multivariate Behrens-Fisher problem.  
Generalizing the univariate 
result described earlier, we shall prove the following result:

\begin{thm}
\label{thm:main}
Suppose that $N_1, N_2 > p$.  Then, almost surely, there are exactly 
$2p+1$ complex solutions of the system of likelihood equations 
for the multivariate Behrens-Fisher problem.  In particular, almost 
surely, there always exists at least one real solution.
\end{thm}


\par


\section{Derivation of the likelihood equations}
\label{section:reduction}
\setcounter{equation}{0} 

Denote by $\bar{X}$ and $\bar{Y}$ the means of the samples from 
$N_p(\mu,\Sigma_1)$ and $N_p(\mu,\Sigma_2)$, respectively.  By standard 
calculations (cf., Mardia, {\it et al.} (1979), p. 142), we find that the 
likelihood equations for estimating $\mu$, $\Sigma_1$ and $\Sigma_2$ are: 
\begin{equation}
\begin{aligned}
\label{MLeqforSigmas}
\Sigmahat_1 & = N_1^{-1}
\sum_{j=1}^{N_1} (X_j-\muhat)(X_j-\muhat)', \\
\Sigmahat_2 & = N_2^{-1}
\sum_{j=1}^{N_1} (Y_j-\muhat)(Y_j-\muhat)'
\end{aligned}
\end{equation}
and
\begin{equation}
\label{MLeqformu}
(N_1\Sigmahat_1^{-1} + N_2\Sigmahat_2^{-1}) \muhat = 
N_1\Sigmahat_1^{-1} \bar{X} + N_2\Sigmahat_2^{-1}\bar{Y}.
\end{equation}

\medskip

\noindent
Some authors have proposed the following iterative algorithm 
for solving (\ref{MLeqforSigmas}) and (\ref{MLeqformu}):

\begin{itemize}

\item[(1)] Begin the iteration with initial estimates 
$\Sigmahat_{i,0} = \wtil{S}_i$, $i=1,2$, where 
\begin{equation}
\begin{aligned}
\label{Stildes}
\wtil{S}_1 & = N_1^{-1}\sum_{j=1}^{N_1} (X_j-\bar{X})(X_j-\bar{X})', \\ 
\wtil{S}_2 & = N_2^{-1}\sum_{j=1}^{N_2} (Y_j-\bar{Y})(Y_j-\bar{Y})'.
\end{aligned}
\end{equation}

\item[(2)] Apply (\ref{MLeqformu}) to calculate $\muhat_0$, the 
corresponding estimate of $\mu$, in the form 
$$
\muhat_0 = (N_1\Sigmahat_{1,0}^{-1} + N_2\Sigmahat_{2,0}^{-1})^{-1}
(N_1\Sigmahat_{1,0}^{-1} \bar{X} + N_2\Sigmahat_{2,0}^{-1}\bar{Y}).
$$

\item[(3)] Use the value of $\muhat_0$ obtained in Step (2) to 
calculate $\Sigmahat_{i,1}$, an updated value of $\Sigmahat_{i,0}$, 
using the formulas 
$$
\Sigmahat_{1,1} = \wtil{S}_1 + (\bar{X}-\muhat_0)(\bar{X}-\muhat_0)', 
\quad
\Sigmahat_{2,1} = \wtil{S}_2 + (\bar{Y}-\muhat_0)(\bar{Y}-\muhat_0)',
$$  
which are a consequence of (\ref{Sigma1hat:eq}) and 
(\ref{Sigma2hat:eq}) below.

\item[(4)] Return to Step (2) and update $\muhat_j$ until the sequences 
$\Sigmahat_{1,j}$ and $\Sigmahat_{2,j}$, $j=1,2,3,\ldots$, converge.
\end{itemize}

We are grateful to Mathias Drton for pointing out that his work in 
Drton and Eichler (2006) implies that this algorithm converges to 
a saddle point or a local (but not necessarily a global) maximum 
of the likelihood function.  
 If the likelihood function were found to be 
multimodal, a phenomenon which has been encountered recently by 
Drton and Richardson (2004) in a study of seemingly unrelated 
regression models, then any numerical algorithm for solving the 
system of likelihood equations necessarily must include some 
information about the choice of initial values.  

At first glance, the likelihood equations may appear to be a 
system of $p(p+2)$ equations in $p(p+2)$ variables comprising 
the $p$ components of $\mu$ and the $p(p+1)/2$ entries of both 
$\Sigma_1$ and $\Sigma_2$.  However, a closer inspection of 
(\ref{MLeqforSigmas}) and (\ref{MLeqformu}) reveals that if 
$\muhat$ is known then $\Sigmahat_1$ and $\Sigmahat_2$ are determined 
completely. We shall show later how to eliminate $\Sigmahat_1$ 
and $\Sigmahat_2$ from (\ref{MLeqformu}) to obtain a system 
of $p$ cubic equations in the variables $\muhat_1,\ldots,\muhat_p$.  

\par 

\begin{prop}
\label{prop:reduction}
The likelihood equations (\ref{MLeqforSigmas}) and (\ref{MLeqformu}) 
for the Behrens-Fisher problem are equivalent to 
\begin{equation}
\label{MLeqforBF}
\frac{N_1 \wtil{S}_1^{-1}(\bar{X}-\muhat)}
{1+(\bar{X}-\muhat)'\wtil{S}_1^{-1}(\bar{X}-\muhat)} + 
\frac{N_2 \wtil{S}_2^{-1}(\bar{Y}-\muhat)}
{1+(\bar{Y}-\muhat)'\wtil{S}_2^{-1}(\bar{Y}-\muhat)} = 0.
\end{equation}
\end{prop}

\begin{proof}
We apply to the sums in (\ref{MLeqforSigmas}) the standard procedure 
of writing each term $X_i-\muhat$ as $X_i-\bar{X} + \bar{X}-\muhat$, 
and similarly for each term $Y_i-\muhat$.  This leads to the formulas  
\begin{equation}
\label{Sigma1hat:eq}
\Sigmahat_1 = \wtil{S}_1 + (\bar{X}-\muhat)(\bar{X}-\muhat)'
\end{equation}
and 
\begin{equation}
\label{Sigma2hat:eq}
\Sigmahat_2 = \wtil{S}_2 + (\bar{Y}-\muhat)(\bar{Y}-\muhat)'
\end{equation}
where $\wtil{S}_1$ and $\wtil{S}_2$ are defined in (\ref{Stildes}).  
By a special case of Woodbury's theorem (cf., Muirhead (1982), p. 580, 
Theorem A5.1) we have, for any nonsingular $p \times p$ matrix $M$ and 
any column vector $v \in \R^p$, 
$$
(M+vv')^{-1} = M^{-1} - \frac{M^{-1}vv'M^{-1}}{1+v'M^{-1}v}.
$$
Multiplying the latter equation on each side from the right by $v$ and 
simplifying the result, we obtain 
\begin{eqnarray*}
(M+vv')^{-1} v & = & M^{-1}v - \frac{M^{-1}vv'M^{-1}v}{1+v'M^{-1}v} \\
& = & \frac{(1+v'M^{-1}v)M^{-1}v - (M^{-1}v)(v'M^{-1}v)}{1+v'M^{-1}v} \\
& = & \frac{M^{-1}v}{1+v'M^{-1}v}.
\end{eqnarray*}
Setting $M = \wtil{S}_1$ and $v = \bar{X}-\muhat$, we obtain 
\begin{eqnarray}
\label{woodbury1}
\Sigmahat_1^{-1}(\bar{X}-\muhat) & \equiv &
\big(\wtil{S}_1 + (\bar{X}-\muhat)(\bar{X}-\muhat)'\big)^{-1} 
(\bar{X}-\muhat) \nonumber \\
& = & \frac{\wtil{S}_1^{-1}(\bar{X}-\muhat)}
{1+(\bar{X}-\muhat)'\wtil{S}_1^{-1}(\bar{X}-\muhat)},
\end{eqnarray}
and, similarly,
\begin{equation}
\label{woodbury2}
\Sigmahat_2^{-1}(\bar{Y}-\muhat) = 
\frac{\wtil{S}_2^{-1}(\bar{Y}-\muhat)}
{1+(\bar{Y}-\muhat)'\wtil{S}_2^{-1}(\bar{Y}-\muhat)}.
\end{equation}
On rewriting (\ref{MLeqformu}) as 
$$
N_1 \Sigmahat_1^{-1}(\bar{X}-\muhat) + 
N_2 \Sigmahat_2^{-1}(\bar{Y}-\muhat) = 0,
$$
it follows from (\ref{woodbury1}) and (\ref{woodbury2}) that 
(\ref{MLeqformu}) is equivalent to (\ref{MLeqforBF}).  
\end{proof}

\section{The maximum likelihood degree of the Behrens-Fisher problem}
\label{MLdegree}
\setcounter{equation}{0}

Following Catanese, {\it et al.} (2006) and Ho\c{s}ten, {\it et al.} 
(2005) we will call the number of complex solutions to the likelihood 
equations the {\em maximum likelihood degree}.  In this section we 
shall prove Theorem \ref{thm:main}, namely, that the maximum 
likelihood (or ML) degree of the Behrens-Fisher problem is $2p+1$.  
Before providing the details of the proof, it is instructive to 
understand why the theorem holds for $p=1$ and $p=2$.  Let us 
denote by $D_X(\muhat)$ and $D_Y(\muhat)$ the denominators 
$1+(\bar{X}-\muhat)'\wtil{S}_1^{-1}(\bar{X}-\muhat)$ and 
$1+(\bar{Y}-\muhat)'\wtil{S}_2^{-1}(\bar{Y}-\muhat)$, respectively, 
which appear in the likelihood equations (\ref{MLeqforBF}).  

\begin{lemma}
\label{lem:bothzero}
Let 
\begin{equation}
\label{BFfinal}
N_1 D_Y(\muhat) \wtil{S}_1^{-1}(\bar{X}-\muhat) +
N_2 D_X(\muhat) \wtil{S}_2^{-1}(\bar{Y}-\muhat) = 0
\end{equation}
be the system of polynomial equations obtained by clearing denominators 
in the likelihood equations (\ref{MLeqforBF}), and suppose that $\muhat$ 
is a solution to (\ref{BFfinal}).  Then $D_X(\muhat) = 0$ if and only if 
$D_Y(\muhat) = 0$.  
\end{lemma}

\begin{proof} Suppose that $D_X(\muhat) = 0$.  On multiplying 
(\ref{BFfinal}) from the left by $(\bar{X}-\muhat)'$ we obtain 
$N_1 D_Y(\muhat)(D_X(\muhat) - 1) = 0$ 
and so we deduce that $D_Y = 0$. Similarly, starting with the 
assumption that $D_Y = 0$, we deduce that $D_X = 0$.  
\end{proof}

We remark that, because $D_X(\muhat)$ and $D_Y(\muhat)$ are strictly 
positive for any real $\muhat$, the system of equations 
(\ref{MLeqforBF}) and (\ref{BFfinal}) are equivalent when determining 
{\em real} solutions 
only.  However, in the calculation of {\em complex} solutions, the 
likelihood equations (\ref{MLeqforBF}) are not equivalent to 
(\ref{BFfinal}) since it is possible that the denominators are 
zero for complex $\muhat$. 

Let $J$ be the ideal defined by the equations (\ref{BFfinal}) and 
let $I = \langle D_X(\muhat),D_Y(\muhat) \rangle$ be the ideal of 
zeros common to the first and second denominators.  Then we need 
to compute and count 
the solutions to $J:I$.  For the case in which $p=1$ there 
is a single univariate cubic polynomial in (\ref{BFfinal}) which, 
generically, has three complex roots.  Since two generic univariate 
polynomials (in this case, $D_X(\muhat)$ and $D_Y(\muhat)$) have no 
common roots then 
the ideal $I$ has, in general, no solutions.  Hence we conclude that 
$J:I$ has exactly three solutions for the case in which $p = 1$.  

We now consider the case in which $p=2$.  Since two quadrics in two 
variables have, generically, four complex roots then there are four 
generic solutions to $I$.  Similarly, since two cubics in two 
variables have generically nine complex roots then there are nine 
generic solutions to $J$.  Therefore $J:I$ has five complex roots 
for the case in which $p = 2$.  

Unfortunately, this counting argument fails even for $p=3$.  In this 
case, we have two quadrics in three variables, so there are infinitely 
many solutions to $I$ and hence also to $J$.  Yet, $J:I$ still has 
finitely many solutions.  Theorem \ref{thm:main} follows from the 
following result of Catanese, {\it et al.} (2006).  

\begin{thm} 
\label{thm1} {\rm (Catanese, {\it et al.}, 2006)}
Let $f_1,\ldots,f_n$ be polynomials of degrees $b_1,\ldots,b_n$, 
respectively, in the variables $x_1,\ldots,x_d$; let 
$u_1,\ldots,u_n$ be integers; let $f = f_1^{u_1} \cdots f_n^{u_n}$; 
and consider the critical equations 
\begin{equation*}
\label{criticalequations}
\frac{1}{f} \frac{\partial f}{\partial x_1} \,\, = \,\,
\frac{1}{f} \frac{\partial f}{\partial x_2} \,\, = \,\,
\cdots \,\, = \,\,
\frac{1}{f} \frac{\partial f}{\partial x_n} \,\, = \,\, 0 
\end{equation*}
of $\log f = \sum_{i=1}^n u_i \log f_i$.  If the number of complex 
solutions to this system of equations is finite then that number is 
less than or equal to the coefficient of $z^d$ in the generating 
function 
\begin{equation*}
\label{thm1formula}
\frac{(1-z)^d }{(1-b_1z) (1-b_2z) \cdots (1-b_nz)}.
\end{equation*}
Equality holds if the coefficients of the polynomials $f_i$ are 
sufficiently generic.
\end{thm}

Before we proceed, there are a few points that need clarification in Theorem 
\ref{thm1}. First of all, given
the integers $b_1, \ldots, b_n$ there exists a fixed polynomial $G= G_{b_1, \ldots, b_n}$ in the coefficients
of $n$ polynomials in $d$ variables with degrees $b_1, \ldots, b_n$ 
so that we call  $f_1, \ldots, f_n$ {\em generic} if $G(f_1, \ldots, f_n) \neq 0$.
Furthermore, when  $f_1, \ldots, f_n$ is generic in this sense
the number of complex solutions to the critical equations is given by the formula in
the statement of the theorem.
In other words, genericity already implies the finiteness of the number
of complex solutions. This follows from Theorem 5 in Catanese, {\it et al.} (2006)
leading up to the proof of Theorem \ref{thm1}.

Before providing the proof of Theorem \ref{thm:main}, we show that 
the coefficients of $D_X(\muhat)$ and $D_Y(\muhat)$ are generic for 
almost all data $X_1,\ldots,X_{N_1}$ and $Y_1,\ldots,Y_{N_2}$.
First we need the following result which has a standard proof in the literature
(for instance based on the argument on page 76 in Anderson (2003)). We present
our own proof.
\begin{lemma}
\label{lem:genericmatrix}
Suppose that $N+1 > p$. Then given any $p \times p$ positive definite 
matrix $S$ there exist $X_1,\ldots,X_{N+1} \in \R^p$ such 
that $S = \sum_{i=1}^{N+1} (X_i - \bar{X})(X_i - \bar{X})'$.
\end{lemma}

\begin{proof} We can assume, without loss of generality that 
$\bar{X} = 0$.  Now let $X_i = (X_{i1}, \ldots, X_{ip})'$ for 
$i = 1, \ldots, N$ and let $X_{N+1} = - \sum_{i=1}^{N} X_i$.  Given 
a positive definite matrix $S$, there exists a nonsingular symmetric 
matrix $U$ such that $USU' = \Lambda$ where $\Lambda$ is a diagonal 
matrix with diagonal entries $\lambda_i > 0$, $i=1,\ldots,p$.  Hence 
it is enough to prove the result for diagonal matrices $\Lambda$.
The required identity $\Lambda = \sum_{i=1}^{N+1} X_i X_i'$ gives 
rise to $p(p+1)/2$ polynomial equations, namely, 
$$
\sum_{i=1}^N \sum_{j=1}^{N} X_{ik}X_{jk} = \frac{\lambda_k}{2},
$$
$k=1,\ldots,p$, and 
$$
X_{1i}(X_{1j} + \sum_{k=1}^N X_{kj}) + X_{2i}(X_{2j} + \sum_{k=1}^N
X_{kj})  + \cdots +
X_{Ni}(X_{Nj} + \sum_{k=1}^N X_{kj})= 0,
$$
$1 \leq i < j \leq p$.  We claim that there exists at least one real 
solution to the above system where $X_{ij} = 0$ for $j=1,\ldots,p$ 
and $i=j+1,\ldots,N$.  It is easy to check that 
$X_{11} = \sqrt{\frac{\lambda_1}{2}}$ with
$X_{1k} = X_{2k} = \cdots = X_{k-1,k} = \sqrt{\frac{\lambda_k}{k(k+1)}}$ 
and $X_{kk} = -k \sqrt{\frac{\lambda_k}{k(k+1)}}$ for $k=2,\ldots,p$ 
give such a solution.
\end{proof}

\begin{thm} 
\label{thm:genericdata}
For generic data $X_1,\ldots,X_{N_1}$ and $Y_1,\ldots,Y_{N_2}$ 
the denominators $D_X$ and $D_Y$ are generic.
\end{thm}

\begin{proof}  $D_X$ and $D_Y$ are quadratic forms in $p$ variables.
In the light of our remarks after Theorem \ref{thm1}
there exists a fixed polynomial $G$ in the coefficients of two quadratic 
forms in $p$ variables such that $D_X$ and $D_Y$ are generic if
$G(D_X, D_Y) \neq 0$.
We need to show that this condition holds for generic data.  As in 
the proof of Lemma \ref{lem:genericmatrix}, the entries of 
$\wtil{S}_1$ and $\wtil{S}_2$ are polynomials in the data.  The same 
lemma implies that the polynomial maps defined from the data spaces 
$\R^{p \times N_1}$ and $\R^{p \times N_2}$ are 
surjective onto the cone of semidefinite matrices in $\R^{p(p+1)/2}$.  
Therefore there exist data vectors $X_1, \ldots, X_{N_1}$ and $Y_1, \ldots, Y_{N_2}$ such that
$G(D_X, D_Y) \neq 0$.
If the statement in the theorem is not true, then there exists a
Zariski open subset $U \subset \R^{p \times N_1} \times \R^{p \times N_2}$ such that
for all $(X_1, \ldots, X_{N_1}\, :\, Y_1, \ldots, Y_{N_2}) \in U$ we have
$G(D_X, D_Y)= 0$. But this means that $G$ is identically
zero, and this is a contradiction.
\end{proof} 

\noindent {\em Proof of Theorem \ref{thm:main}:}  Denoting by 
$L(\mu,\Sigma_1,\Sigma_2)$ the likelihood function for the 
Behrens-Fisher problem, then it is well-known that 
$$
L(\muhat,\Sigmahat_1,\Sigmahat_2) = (2\pi e)^{-(N_1+N_2)p/2} \, 
|\Sigmahat_1|^{-N_1/2} \, |\Sigmahat_2|^{-N_2/2}.
$$
By (\ref{Sigma1hat:eq}) and (\ref{Sigma2hat:eq}), we have 
$|\Sigmahat_1| = |\wtil{S}_1| \cdot  D_X(\muhat)$ and 
$|\Sigmahat_2| = |\wtil{S}_2| \cdot  D_Y(\muhat)$.  Therefore
\begin{multline*}
L(\muhat,\Sigmahat_1,\Sigmahat_2) \\ = (2\pi e)^{-(N_1+N_2)p/2} \, 
|\wtil{S}_1|^{-N_1/2} \, |\wtil{S}_2|^{-N_2/2} 
\big(D_X(\muhat)\big)^{-N_1/2} \, \big(D_Y(\muhat)\big)^{-N_2/2}.
\end{multline*}
It now is clear that, to find the maximum value of $L$, we need 
to minimize 
\begin{equation}
\label{Lmaximum}
\big(1 + (\bar{X}-\muhat)'\wtil{S}_1^{-1}(\bar{X}-\muhat)\big)^{N_1/2} \, 
\big(1 + (\bar{Y}-\muhat)'\wtil{S}_2^{-1}(\bar{Y}-\muhat)\big)^{N_2/2}.
\end{equation}
Equivalently, we may minimize the logarithm of this expression, and 
since the critical equations of the logarithm of (\ref{Lmaximum}) 
are precisely the likelihood equations in Proposition 
\ref{prop:reduction}, then Theorem \ref{thm1} implies that the 
maximum likelihood degree of the Behrens-Fisher problem is equal 
to the coefficient of $z^p$ in the power series expansion of 
the rational function
$$
\frac{(1-z)^p}{(1-2z)^2},
$$ 
provided that the data $X_1, \ldots, X_{N_1}$ and
$Y_1, \ldots, Y_{N_2}$, and hence $D_X$ and $D_Y$ are generic.  By expanding this rational 
function in a power series in $z$, we find that this coefficient 
equals 
$$
\sum_{i+j = p} (-1)^i  \, 2^j \, {p \choose i} (j+1),
$$
and an elementary calculation shows that this sum equals $2p+1$. 

Step 4 in Algorithm 7 in Ho\c{s}ten, {\it et al.} 
(2005) and the theory of Gr\"obner bases imply that 
all $2p+1$ complex solutions can be obtained from the roots of a 
univariate polynomial of degree $2p+1$. Since $D_X(\muhat)$ and 
$D_Y(\muhat)$ have real 
coefficients, this univariate polynomial also has real coefficients.  
In particular, since roots occur in complex conjugate pairs then 
at least one root is real. \hfill $\Box$

\begin{rmk}{\rm 
We note that our arguments which led to the derivation of the 
ML degree of the Behrens-Fisher problem also apply to the more 
general problem of multivariate analysis of variance (MANOVA).  
Suppose that we have independent multivariate normal populations 
$N_p(\mu,\Sigma_1),\ldots,N_p(\mu,\Sigma_{k+1})$ and that, 
on the basis of random samples from each population, we wish to 
derive the maximum likelihood estimators of the parameters $\mu$ 
and $\Sigma_1,\ldots,\Sigma_{k+1}$.  
By arguments similar to those in Section \ref{section:reduction}, 
we obtain analogous likelihood equations as in Proposition 
\ref{prop:reduction} where now there are $k+1$ rational summands 
in each of the $p$ equations.  It then follows from Theorem 
\ref{thm1} that the ML degree for the MANOVA problem is 
\begin{equation}
\label{MLdegMANOVA}
d(k,p) := \sum_{i+j=p} (-1)^i \, 2^j \, {p \choose i} \, {j+k \choose k}.
\end{equation}
By writing this result in the form 
$$
d(k,p) = 1 + \sum_{j=1}^p (-1)^{p-j} \, 2^j \, {p \choose j} \, 
{j+k \choose k},
$$
we find that $d(k,p)$ is odd; therefore, there always exists 
a real solution to the system of likelihood equations. 

We note that $d(k,p)$ can be evaluated using methods from the 
calculation of combinatorial sums, as follows:  First, we write 
$$
2^j \, {j+k \choose k} = \frac{1}{k!} \, \Big(\frac{d}{dt}\Big)^k \, 
t^{j+k} \,\Bigg|_{t=2} \ .
$$
Inserting this formula in the sum in (\ref{MLdegMANOVA}) and 
interchanging derivatives and summation, we obtain 
\begin{eqnarray}
\label{MANOVAdeg}
d(k,p) & = & \frac{1}{k!} \, \Big(\frac{d}{dt}\Big)^k \, t^k \, 
\sum_{j=0}^p \, (-1)^{p-j} \, {p \choose j} \, t^j \,\Bigg|_{t=2} 
\nonumber \\
& = & \frac{1}{k!} \, \Big(\frac{d}{dt}\Big)^k \, t^k \, (t-1)^p 
\,\Big|_{t=2} \ .
\end{eqnarray}
In particular, $d(1,p) = 2p+1$, the ML degree of the Behrens-Fisher 
problem, and $d(2,p) = 2p(p+1)+1$.  The general formula for $d(k,p)$ 
is interesting even in the case $p = 1$, for it yields the ML degree 
of the one-dimensional $(k+1)$-population MANOVA problem to be 
$2k+1$.  Further, by substituting $t=(1+u)/2$ in (\ref{MANOVAdeg}), 
we recognize the outcome as Rodrigues' formula (Szeg\"o, 1939, p. 66) 
for a Jacobi polynomial $P^{(p-k,0)}_k$, and we obtain 
$d(k,p) = P^{(p-k,0)}_k(3)$, $p \ge k$.}
\end{rmk}

\section{Simulations and a large sample size result}
\label{simulations}
\setcounter{equation}{0}

Having determined the number of solutions of the system of 
likelihood equations (\ref{BFfinal}) it is natural to seek the 
number of {\it real} solutions, for it is those solutions which 
are of interest in statistical inference.  Not surprisingly, it 
appears to be difficult to determine an algebraic 
expression for the number of real solutions of the system; 
indeed, this is also the case for the general theory of 
systems of polynomial equations.

To study the real solutions of the system (\ref{BFfinal}), we 
considered the case in which $p = 2$, presenting empirical evidence 
that multiple solutions occur rarely if the model is correctly 
specified.  In each simulation run, we 
first used a random number generator to generate sample sizes 
$N_1$ and $N_2$, and a mean vector $\mu$.  We next generated lower 
triangular matrices $T_1$ and $T_2$ with positive diagonal entries, 
after which we set $\Sigma_k = T_kT_k'$, $k = 1, 2$.  Finally, we 
simulated a random sample of vectors $Z_1,\ldots,Z_{N_1}$ from 
$N_2(0,I_2)$, and then we set $X_j = T_1Z_j + \mu$, 
$j=1,\ldots,N_1$.  It follows from standard distribution theory 
that $X_1,\ldots,X_{N_1}$ constitutes a simulated sample from the 
bivariate normal population $N_2(\mu,\Sigma_1)$.  In a similar 
manner, we simulated an independent random sample 
$Y_1,\ldots,Y_{N_1}$ from $N_2(\mu,\Sigma_2)$.  

The solutions of the resulting likelihood 
equations (\ref{BFfinal}) were computed numerically using 
{\tt{PHCpack}} (Verschelde, 1999), a software package which 
implements polyhedral homotopy continuation methods for solving 
systems of polynomial equations.  
The results of our simulations show that multiple solutions can 
occur.  For example, for $N_1=11$, $N_2=5$, and the summary 
statistics 
$$
\begin{array}{lclrcl}
\bar{X} & = & 
\left(
\begin{array}{c}
-1.5516 \\ -9.4713 
\end{array} 
\right), 
\qquad
\wtil{S}_1 & = & 
\left(
\begin{array}{rr}
0.3998  &   -0.1026 \\ -0.1026 & 0.2378 
\end{array} 
\right), 
\\ \\
\bar{Y} & = & 
\left(
\begin{array}{r}
-1.9175 \\ -10.4805 
\end{array} 
\right), 
\qquad
\wtil{S}_2 & = & 
\left(
\begin{array}{rr}
0.4193  &    0.0792 \\  0.0792 & 0.0334 
\end{array} 
\right),
\end{array}
$$
the real solutions for $\mu$ are 
$$
\left(
\begin{array}{c}
-1.3570 \\ -10.2957
\end{array}
\right), \quad 
\left(
\begin{array}{c}
-1.2478 \\ -9.9902
\end{array}
\right), \quad \hbox{and} \quad
\left(
\begin{array}{c}
-1.4451 \\ -9.6333
\end{array}
\right).
$$
This example seems, however, to be a rare exception.  Indeed, 
we found that the bivariate Behrens-Fisher likelihood 
equations (\ref{BFfinal}) had one real solution in about 
99.5\% of simulations, three real solutions in about 0.5\% of 
simulations, and we found no instances in which the equations had 
five real solutions.  However, it is possible that (\ref{BFfinal}) has
five real solutions when the data is generated from a "wild" distribution
and not from the corresponding multivariate distributions. 
For instance,  for $N_1=15$, $N_2=28$, and 
$$
\begin{array}{lclrcl}
\bar{X} & = & 
\left(
\begin{array}{c}
-4 \\ -3
\end{array} 
\right), 
\qquad
\wtil{S}_1 & = & 
\left(
\begin{array}{rr}
49.3619  &   -45.0547 \\ -45.0547 & 42.4495 
\end{array} 
\right), 
\\ \\
\bar{Y} & = & 
\left(
\begin{array}{r}
4 \\ 1 
\end{array} 
\right), 
\qquad
\wtil{S}_2 & = & 
\left(
\begin{array}{rr}
52.8534  &    19.8380 \\  19.8380 & 9.0472
\end{array} 
\right),
\end{array}
$$
the real solutions for $\mu$ are 
$$
\left(
\begin{array}{c}
3.9822 \\ 1.0443
\end{array}
\right), \quad 
\left(
\begin{array}{r}
-3.7286 \\ 3.2906
\end{array}
\right), \quad 
\left(
\begin{array}{r}
-2.4192 \\ 4.6925
\end{array}
\right),
\quad 
\left(
\begin{array}{c}
2.0437 \\ 5.8993
\end{array}
\right), \quad 
\left(
\begin{array}{c}
1.0089 \\ 8.2001
\end{array}
\right).
$$

\vskip 0.3cm
To test for distinctions between the case of small and large samples 
in the bivariate case, we performed simulations in which $N_1$ and 
$N_2$ were randomly generated (uniform distribution) between 3 and 15.  The outcomes are 
given as follows, with percentages rounded-off to two decimal places:

\begin{center}
\vskip 5pt
\framebox{Please place Table 1 here}
\vskip 12pt
\end{center}

\noindent
As noted above, none of these simulation resulted in five real 
solutions.  

In the case of larger samples, our simulations resulted in the 
following outcomes: 

\begin{center}
\vskip 5pt
\framebox{Please place Table 2 here}
\vskip 12pt
\end{center}

\noindent
Here again, no simulation resulted in five real solutions.  (In both cases, the population mean $\mu$ is randomly generated from a
uniform distribution on the subspace $[-20,20] \times [-20,20]$, and the
population covariance matrices $\Sigma_1$ and $\Sigma_2$ are randomly
generated in the manner described above, with positive diagonal entries
whose values are no greater than 10.)

In summary, there seems to be little chance that a randomly generated, 
two-dimensional Behrens-Fisher problem will have three or more real 
solutions, and there is a high chance that it will have a unique 
real solution.  The following supports the second conclusion for
large sample sizes.



\begin{thm} 
\label{thm:uniqueness} 
Suppose that the random samples $X_1,\ldots,X_{N_1}$ and 
$Y_1,\ldots,Y_{N_2}$ are drawn from independent normal populations 
$N_p(\mu,\Sigma_1)$ and $N_p(\mu,\Sigma_2)$, respectively.
As $N_1, N_2 \to \infty$ the likelihood equations (\ref{MLeqforBF})
for the Behrens-Fisher problem has a unique real root with probability one.
\end{thm}
\begin{proof}  If $\bar{X} = \bar{Y}$ then 
it follows from (\ref{Lmaximum}) that the {\it unique} real solution 
of the likelihood equations is $\muhat = \bar{X} = \bar{Y}$.  
Without loss of generality we can assume that $\bar{X} = \bar{Y} = 0$,
and with this the likelihood equations are
\begin{equation}
\label{MLeqforunique}
\frac{N_1 \wtil{S}_X^{-1}\muhat}
{1+\muhat'\wtil{S}_X^{-1}\muhat} + 
\frac{N_2 \wtil{S}_Y^{-1}\muhat}
{1+\muhat'\wtil{S}_Y^{-1}\muhat} = 0.
\end{equation}
We argue that $\muhat = 0$ is a solution of multiplicity one for the system 
obtained by clearing denominators in (\ref{MLeqforunique}). Let $I$ be the ideal
in $\mathbb{C}[\mu_1, \ldots, \mu_p]$ generated by these $p$ equations. 
The multiplicity of $\muhat = 0$ is the length of the artinian module
$$\frac{\mathbb{C}[\mu_1, \ldots, \mu_p]_{\langle \mu_1, \ldots, \mu_p \rangle}}{I \cdot \mathbb{C}[\mu_1, \ldots, \mu_p]_{\langle \mu_1, \ldots, \mu_p \rangle}}$$
over the local ring $\mathbb{C}[\mu_1, \ldots, \mu_p]_{\langle \mu_1, \ldots, \mu_p \rangle}$.
$I$ is generated by $p$ polynomials given by 
$(N_1\wtil{S}_X^{-1} + N_2\wtil{S}_Y^{-1}) \muhat + N_1(\muhat'\wtil{S}_Y^{-1}\muhat) \wtil{S}_X^{-1} \muhat +  N_2(\muhat'\wtil{S}_X^{-1}\muhat) \wtil{S}_Y^{-1} \muhat = 0 $.
Each of these polynomials consists of a linear term and a cubic term. 
With probability one the rank of $N_1\wtil{S}_X^{-1} + N_2\wtil{S}_Y^{-1}$ over 
$\mathbb{C}$ is $p$, 
and hence we can assume that $I$ is generated by $p$ polynomials 
of the form $ \mu_i + g_i$ where $g_i$ has degree three.
This implies that the initial ideal of $I$ in the local ring $\mathbb{C}[\mu_1, \ldots, \mu_p]_{\langle \mu_1, \ldots, \mu_p \rangle}$ with respect to the local term order 
{\em anti-graded revlex} as on page 152 of Cox, {\it et al.} (1998) is
$\langle \mu_1, \ldots, \mu_p \rangle$. By Corollary 4.5 of Cox, {\it et al.} (1998), we conclude
that the length of the above module and hence the multiplicity of $\muhat = 0 $ is one.
Now as $N_1, N_2 \to \infty$, by  the Law of Large Numbers, $\bar{X}$ and $\bar{Y}$ 
converge to $\mu$, and $S_1$ and $S_2$ converge to $S_X$ and $S_Y$.  Since $\muhat = 0$
is the unique real solution to   (\ref{MLeqforunique}) with multiplicity one, and by the 
continuity of solutions to the general likelihood equations (\ref{MLeqforBF}), we conclude
that with probability one (\ref{MLeqforBF}) has a unique real solution.
\end{proof}



\vskip 0.5cm

\noindent {\large\bf Acknowledgment}.  We thank Mathias Drton, Bernd 
Sturmfels, and the Editors for discussions and comments on initial 
drafts of this manuscript.  Richards' work was supported in part by 
NSF grant DMS-0705210.

\par

\bigskip
\bigskip
\bigskip

\noindent{\large\bf References}
\begin{description}

\item
Anderson, T. W. (2003).  {\it An Introduction to Multivariate 
Statistical Analysis}, third edition.  Wiley, New York.

\item
Behrens, W. U. (1929).  Ein Beitrag zur Fehlerberechnung bei wenigen 
Beobachtungen.  {\it Landwirtschaftliche Jahrb\"ucher}, {\bf 68}, 
807--837.


\item
Buot, M.-L. G. and Richards, D. St. P. (2006).  Counting and locating 
the solutions of polynomial systems of maximum likelihood equations, I. 
{\it J. Symbolic Comput.} {\bf 41}, 234--244.

\item Catanese, F., Ho\c{s}ten, S., Khetan, A., and Sturmfels, B. 
(2006).  The maximum likelihood degree.  {\it Amer. J. Math.} 
{\bf 128}, 671--697.

\item
Cox, D. A., Little, J., and O'Shea, D. (1998).  
{\it Using Algebraic Geometry}.  Springer, New York.

\item
Drton, M. (2007).  Multiple solutions to the likelihood equations 
in the Behrens-Fisher problem. {\tt http://arxiv.org/abs/0705.4516.}

\item
Drton, M. and Eichler, M. (2006). Maximum likelihood estimation in 
Gaussian chain graph models under the alternative Markov property. 
{\it Scandinavian Journal of Statistics} {\bf 33}, 247--257.

\item
Drton, M. and Richardson, T. (2004).  Multimodality of the 
likelihood in the bivariate seemingly unrelated regressions 
model.  {\it Biometrika} {\bf 91}, 383--392.

\item
Fisher, R. A. (1939).  The comparison of samples with possibly 
unequal variances.  {\it Ann. Eugen.}, {\bf 9}, 174--180.

\item 
Ho\c{s}ten, S., Khetan, A., and Sturmfels, B. (2005). Solving 
the likelihood equations. {\it Found. Comput. Math.} {\bf 5}, 
389--407.



\item
Kim, S.-H. and Cohen, A. (1998).  On the Behrens-Fisher problem: a 
review.  {\it J. Educational Behavioral Statist.}, {\bf 23}, 356-377.


\item
Linnik, Yu. V. (1967).  On the elimination of nuisance parameters in 
statistical problems.  {\it Proc. Fifth Berkeley Sympos. Math. Statist. 
and Probability}, Vol. I: Statistics, pp. 267--280, Univ. California 
Press, Berkeley, CA.

\item
Mardia, K. V., Kent, J. T., and Bibby, J. M. (1979).  
{\it Multivariate Analysis}.  Academic Press, New York.

\item
Pachter, L. and Sturmfels, B. (2005).  {\it Algebraic Statistics 
for Computational Biology}.  Cambridge University Press, 
Cambridge, U.K.

\item
Stuart, A. and Ord, J. K. (1994).  {\it Kendall's Advanced Theory 
of Statistics}, 6th edition.  Edward Arnold, London.

\item
Sugiura, N. and Gupta, A. K. (1987).  Maximum likelihood estimates 
for the Behrens-Fisher problem.  {\it J. Japan Statist. Soc.}, 
{\bf 17}, 55--60.

\item
Szeg\"o, G. (1939).  {\it Orthogonal Polynomials}.  American 
Mathematical Society, New York, NY.

\item
Verschelde, J. (1999).  Algorithm 795: PHCpack: A general-purpose 
solver for polynomial systems by homotopy continuation. 
{\it ACM Trans.  Math. Software}, {\bf 25}, 251--276.

\item
Wallace, D. L. (1980).  The Behrens-Fisher and Fieller-Creasy 
problems.  {\it R. A. Fisher: An Appreciation}, pp. 119--147.  
Lecture Notes in Statist., {\bf 1}, Springer, New York.

\end{description}

\newpage

\begin{table}[!h]
\caption{Simulations with $3 \le N_1, N_2 \le 15$} 
\begin{center}
\renewcommand{\arraystretch}{1.0}
\begin{tabular}{|c|r|r|}
\hline
Number of solutions & Frequency & Percentage \\
\hline
1 & 4450 & 99.29\% \\
3 &   32 &  0.71\% \\
\hline
\end{tabular}
\end{center}
\end{table}

\vskip 1truein

\begin{table}[!h]
\caption{Simulations with $15 \le N_1, N_2 \le 60$} 
\begin{center}
\renewcommand{\arraystretch}{1.1}
\begin{tabular}{|c|r|r|}
\hline
Number of solutions & Frequency & Percentage \\
\hline
1 & 4404 & 99.46\% \\
3 &   24 &  0.54\% \\
\hline
\end{tabular}
\end{center}
\end{table}

\newpage

\vskip .65cm
\noindent
Department of Mathematics and Computer Science, Xavier University, 
Cincinnati, OH 45207.
\vskip 2pt
\noindent
E-mail: buotm@xavier.edu
\vskip 10pt
\noindent
Department of Mathematics, San Francisco State University,
1600 Holloway Avenue, San Francisco, CA 94132. 
\vskip 2pt
\noindent
E-mail: serkan@math.sfsu.edu
\vskip 10pt
\noindent
Department of Statistics, Penn State University, University Park, 
PA 16802, and the Statistical and Applied Mathematical Sciences 
Institute, Research Triangle Park, NC 27709.
\vskip 2pt
\noindent
E-mail: richards@stat.psu.edu

\end{document}